\documentstyle{article}
\begin{document}

\newcommand{\HDS}{\vrule width0pt height2.3ex depth1.05ex\displaystyle}

\def\f#1#2{{{\HDS #1}\over{\HDS #2}}}

\def\lpravilo#1{ \makebox[-.5em][r]{\mbox{\it #1}} {\mbox{\hspace{0.5em}}}}

\def\str{\mbox{\boldmath$\rightarrow$}}

\newcommand{\qed}{\hfill $\Box$}

\def\k#1#2{\stackrel{\raisebox{-2pt}{\mbox{\tiny $#1$}}}{k}^
{\raisebox{-8pt}{\scriptsize #2}}}

\def\D{\mbox{$\cal D$}}

\def\P{\mbox{$\cal P$}}

\def\L{\mbox{$\cal L$}}

\def\M{\mbox{$\cal M$}}

\def\F{\mbox{$\cal F$}}

\def\W{\mbox{$W$}}

\def\C{\mbox{$\cal C$}}

\def\sp{\mbox{\it SplPre}}

\def\G{\mbox{\it Gen}}

\def\R{\mbox{\it Rel}}

\def\mj{\mbox{$\boldmath 1$}}

\def\ref{\mbox{\rm Ref}}

\def\tr{\mbox{\rm Tr}}

\def\cl{\mbox{\rm Cl}}

\title{\bf A Brauerian Representation\\ of Split Preorders}

\author{{\sc Kosta Do\v sen} and {\sc Zoran Petri\' c}
\\[.05cm]
\\Mathematical Institute, SANU \\
Knez Mihailova 35, p.f. 367 \\
11001 Belgrade, Yugoslavia \\
email: \{kosta, zpetric\}@mi.sanu.ac.yu}
\date{}
\maketitle

\begin{abstract}
\noindent Split preorders are preordering relations on a domain whose
composition is defined in a particular way by splitting the domain into
two disjoint subsets. These relations and the associated composition
arise in categorial proof theory in connection with coherence theorems.
Here split preorders are represented isomorphically
in the category whose arrows are
binary relations, where composition is defined in the usual way.
This representation is related to a classical result of representation
theory due to Richard Brauer.

\vspace{0.3cm}

\end{abstract}

\section{Introduction}

A split preorder is a preordering, i.e. reflexive and transitive, relation
$R\subseteq \W^2$ such that \W\ is equal to the disjoint union of $X$
and $Y$. Every preorder may be conceived as a split preorder, but split
preorders are not composed in the ordinary manner. The set $X$ is
conceived as a domain, $Y$ as a codomain, and composition of split
preorders is defined in a manner that takes this into account. The
formal definition of this composition is not quite straightforward, and
we will not give it before the next section. It is, however, supported
by geometric intuitions (see the end of the next section).

Split preorders and their compositions are interesting for logic because they
arise as relations associated to graphs that have been used in
categorial proof theory for coherence results, which settled questions
of identity criteria for proofs (see \cite{D02} for a survey of this
topic). We will consider this matter in Section 5.

A particular brand of split preorder is made of split equivalences,
namely those split preorders that are equivalence relations, which in
categorial proof theory represent {\it generality} of proofs. The paper
\cite{DP02a} is devoted to this matter.

It was shown in \cite{DP02a} that the category whose arrows are split
equivalences on finite ordinals can be represented isomorphically in
the category whose arrows are binary relations between finite ordinals,
where composition is defined in the usual simple way. This
representation is related to Brauer's representation of Brauer algebras
\cite{B37}, which it generalizes in a certain sense (see \cite{DP02a},
Section 6).

In this paper we will generalize the results of \cite{DP02a}. We will
represent the category whose arrows are split preorders in general
in the category
\R\ whose arrows are binary relations between certain sets. If our
split preorders are relations on finite sets \W, we may, as in
\cite{DP02a}, take
that in \R\ we have as arrows binary relations between finite ordinals. Our
representation in \R, which generalizes \cite{DP02a}, generalizes
further Brauer's representation, and hence its name.

In the next section we will introduce the category \sp\ whose arrows
are split preorders. Section 3 and 4 are devoted to representing
isomorphically \sp\ in \R, and also provide a proof that \sp\ is indeed
a category. Section 3 is about some general matters concerning the
representation of arbitrary preordering relations, which we apply in
Section 4. Section 5 is about matters of categorial proof theory, which
we mentioned above. We give some examples of deductive systems
covering fragments of propositional logic, and of the split preorders
associated with proofs.

\section{The category \sp}

Let \M\ be a family of sets, for whose members we use $X,Y,Z,\ldots$
For $i\in\{s,t\}$ (where $s$ stands for ``source'' and $t$ for
``target''), let $\M^i$ be a family of sets in one-to-one
correspondence with \M. We denote by $X^i$ the element of $\M^i$
corresponding to $X\in\M$. We assume further that for every  $X\in \M$
there is a bijection $i_X:X\str X^i$. Finally, we assume that for
every $U\in\M$, $V\in\M^s$ and $W\in\M^t$, the sets $U$, $V$ and $W$
are mutually disjoint.

For $X,Y\in\M$, let a {\it split relation} of \M\ be a triple $\langle
R,X,Y\rangle$ such that $R\subseteq (X^s\cup Y^t)^2$. The set $X^s\cup
Y^t$ may be conceived as the disjoint union of $X$ and $Y$. A split
relation $\langle R,X,Y\rangle$ is a {\it split preorder} iff $R$ is a
preorder, i.e. a reflexive and transitive relation. As usual, we write
sometimes $xRy$ for $(x,y)\in R$.

We will now build a category called \sp, whose set of objects will be
\M, and whose arrows will be split preorders of \M. For such a split
preorder $\langle R,X,Y\rangle$ we take that $X$ is the source, $Y$ the
target, and we write, as usual, $R:X\str Y$.

The identity arrow $\mj_X:X\str X$ of \sp\ is the split preorder that for
every $i,j\in\{ s,t\}$ and every $u\in X^i$ and every $v\in X^j$
satisfies
\[
(u,v)\in \mj_X \quad {\mbox{\rm iff}} \quad i^{-1}_X(u)=j^{-1}_X(v).
\]

To define composition of arrows in \sp\ is a more involved matter, and
for that we need some auxilliary notions. For every $X,Y\in\M$, let the
function $\varphi^s$ from $X\cup Y^t$ to $X^s\cup Y^t$ be defined by
\[
\varphi^s(x)=\left\{\begin{array}{lll}
s_X(x) & {\mbox{\rm if}} &  x\in X
\\
x & {\mbox{\rm if}} & x\in Y^t,
\end{array}\right.
\]
and let the function $\varphi^t$ from $X^s\cup Y$ to $X^s\cup Y^t$ be
defined by
\[
\varphi^t(x)=\left\{\begin{array}{lll}
x & {\mbox{\rm if}} &  x\in X^s
\\
t_Y(x) & {\mbox{\rm if}} & x\in Y.
\end{array}\right.
\]

For a split relation $R:X\str Y$, let the relations $R^{-s}\subseteq(X\cup
Y^t)^2$ and $R^{-t}\subseteq (X^s\cup Y)^2$ be defined by
\[(x,y)\in R^{-i} \quad {\mbox{\rm iff}} \quad (\varphi^i(x),\varphi^i(y))\in
R
\]
for $i\in\{s,t\}$.
Finally, for an arbitrary binary relation $R$, let $\tr(R)$ be
the transitive closure of $R$.

Then for split preorders $R:X\str Y$ and $P:Y\str Z$ of \M\ we define
their composition $P\ast R:X\str Z$ by
\[
P\ast R=_{\mbox{\scriptsize\it def}}
\tr (R^{-t}\cup P^{-s})\cap (X^s\cup Z^t)^2.
\]
It is clear that $P\ast R:X\str Z$ is a split preorder of \M.

For \sp\ to be a category we need that for $R:X\str Y$ the equations
$R\ast \mj_X=\mj_Y \ast R=R$ hold, and that $\ast$ is associative, which it
is rather complicated to check directly. We will not try to do that
here. So, for the time being, we don't know yet whether \sp\ is a
category. We know only that it is a graph (a set of objects and a set
of arrows) with a family of arrows $\mj_X$ for every object $X$, and with
a binary partial operation on arrows $\ast$. This kind of structure is
called a {\it deductive system}, according to \cite{LS86}. We prove
that \sp\ is a category in the sections that follow.

The {\it strictification} of a preorder $R\subseteq \W^2$ is the
relation $R'\subseteq \W^2$ such that $x R' y$ iff $x R y$ and $x\neq y$.
The strictification of a preorder is an irreflexive relation that
satisfies {\it strict transitivity}:
\[
(xR' y\; \& \; yR' z \; \& \; x\neq z)\Rightarrow xR' z.
\]
If $\ref(R)$ is the reflexive closure of $R$, it is clear that for
a preorder $R$ we have $\ref(R')=R$.
Conversely, for every irreflexive strictly transitive
relation $P \subseteq \W^2$
we have that $\ref(P) \subseteq \W^2$ is a preorder, whose strictification
$(\ref(P))'$ is equal to $P$.

So we may represent preorders by their strictifications. This we do
when we draw diagrams representing split preorders. The composition
$\ast$ of two split preorders is illustrated in the following diagrams:

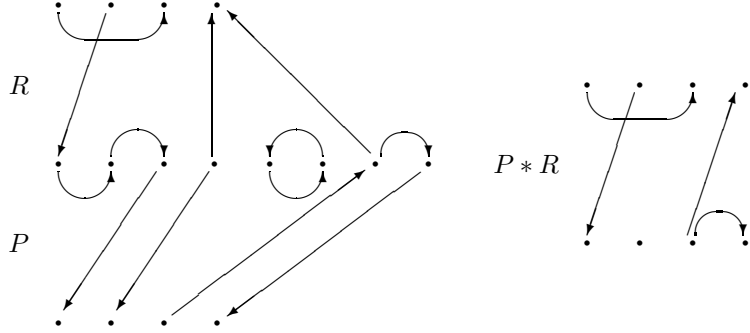
\begin{center}
\begin{picture}(280,140)

\put(57,67){\vector(-2,-3){35}}
\put(77,67){\vector(-2,-3){35}}
\put(63,13){\vector(4,3){73}}
\put(157,67){\vector(-4,-3){73}}
\put(38,127){\vector(-1,-3){18}}
\put(78,73){\vector(0,1){54}}

\put(138,74){\vector(-1,1){54}}

\put(40,65){\vector(0,1){2}}
\put(60,125){\vector(0,1){2}}
\put(120,65){\vector(0,1){2}}
\put(260,95){\vector(0,1){2}}
\put(60,75){\vector(0,-1){2}}
\put(100,75){\vector(0,-1){2}}
\put(160,75){\vector(0,-1){2}}
\put(279,45){\vector(0,-1){2}}

\put(238,97){\vector(-1,-3){18}}
\put(258,43){\vector(1,3){18}}

\put(20,10){\circle*{2}}
\put(40,10){\circle*{2}}
\put(60,10){\circle*{2}}
\put(80,10){\circle*{2}}
\put(20,70){\circle*{2}}
\put(40,70){\circle*{2}}
\put(60,70){\circle*{2}}
\put(79,70){\circle*{2}}
\put(100,70){\circle*{2}}
\put(120,70){\circle*{2}}
\put(140,70){\circle*{2}}
\put(160,70){\circle*{2}}
\put(20,130){\circle*{2}}
\put(40,130){\circle*{2}}
\put(60,130){\circle*{2}}
\put(80,130){\circle*{2}}

\put(220,40){\circle*{2}}
\put(240,40){\circle*{2}}
\put(260,40){\circle*{2}}
\put(280,40){\circle*{2}}
\put(220,100){\circle*{2}}
\put(240,100){\circle*{2}}
\put(260,100){\circle*{2}}
\put(280,100){\circle*{2}}

\put(30,67){\oval(20,20)[b]}
\put(50,73){\oval(20,20)[t]}
\put(110,67){\oval(20,20)[b]}
\put(110,73){\oval(20,20)[t]}
\put(151,73){\oval(18,18)[t]}
\put(40,127){\oval(40,20)[b]}

\put(240,97){\oval(40,20)[b]}
\put(270,43){\oval(18,18)[t]}

\put(5,40){\makebox(0,0){$P$}}
\put(5,100){\makebox(0,0){$R$}}
\put(197,70){\makebox(0,0){$P\ast R$}}

\end{picture}
\end{center}

Every binary relation between $X^s$ and $Y^t$ may be viewed as the
strictification of a split preorder. The composition of such split preorders,
which correspond to binary relations between members of \M, is then
a simple matter: it corresponds exactly to composition of binary relations.
If in \sp\ we keep as arrows only these split preorders, we
obtain a category isomorphic to the category
whose arrows are binary relations between members of \M.
Problems with composition of split
preorders arise if they don't correspond to binary relations between
members of \M.

\section{Representing preorders by sets of functions}

In this section we will consider some general matters concerning the
representation of arbitrary preordering relations. This will serve for
demonstrating in the next section the isomorphism of our representation
of \sp\ in the category whose arrows are binary relations.

Let \W\ be an arbitrary set, and let $R\subseteq\W^2$. Let $p$ be a set
such that the ordinal $2=\{0,1\}$ is a subset of $p$, and let $\leq$ be
a binary relation on $p$ which restricted to 2 is the usual ordering of
2. For $x,y\in\W$, let the function $f_x:\W\str p$ be defined as follows:
\[
f_x(y)=\left\{\begin{array}{lll}
1 & {\mbox{\rm if}} & x R y
\\
0 & {\mbox{\rm if not}} & x R y.
\end{array}
\right.
\]
The function $f_x$ is the characteristic function of the $R$-cone over
$x$. Consider also the following set of functions:
\[
\F(R)=_{\mbox{\scriptsize\it def}} \{f:\W\str p\mid (\forall
x,y\in\W)(xRy\Rightarrow f(x)\leq f(y))\}.\]
We can then establish the following propositions.

\vspace{0.3cm}

\noindent {\sc Proposition} 1.\quad {\it The relation $R$ is reflexive iff
$(\forall x\in \W)f_x(x)=1$.}

\vspace{0.3cm}

\noindent {\sc Proposition} 2.\quad {\it The relation $R$ is transitive iff
$(\forall x\in \W) f_x\in \F(R)$.}

\vspace{0.2cm}

\noindent{\it Proof}. \quad $(\Rightarrow)$ Suppose $yRz$. If $f_x(y)=0$, then
$f_x(y)\leq f_x(z)$, and if $f_x(y)=1$, then $f_x(z)=1$ by the
transitivity of $R$.

\vspace{0.15cm}

\noindent $(\Leftarrow)$ If $yRz\Rightarrow f_x(y)\leq f_x(z)$, then
$yRz\Rightarrow (f_x(y)=0\;{\mbox{\it or }}f_x(z)=1)$, which means
$yRz\Rightarrow(xRy\Rightarrow xRz)$.
\qed

\vspace{0.3cm}

\noindent {\sc Proposition} 3. \quad {\it If $R$ is a preorder, then}
\[
(\ast)\quad (\forall x,y\in\W)(xRy \Leftrightarrow (\forall f\in\F(R))
f(x)\leq f(y)).
\]

\vspace{0.2cm}

\noindent {\it Proof}.\quad Note first that in $(\ast)$ the left-to-right
implication is satisfied by definition. If $(\forall f\in\F(R))f(x)\leq
f(y)$, then $1\leq f_x(y)$ by the left-to-right directions of
Propositions 1 and 2, and so $xRy$.
\qed

\vspace{0.3cm}

\noindent {\sc Proposition} 4.\quad {\it If $\leq$ is a preorder, then $R$ is a
preorder iff $(\ast)$.}

\vspace{0.2cm}

\noindent {\it Proof}.\quad Suppose $\leq$ is a preorder. Then we obtain the
reflexivity of $R$ by taking $x=y$ in $(\ast)$ and by using the
reflexivity of $\leq$.

For the transitivity of $R$, suppose $xRy$ and $yRz$. Then for every
$f\in\F(R)$ we have $f(x)\leq f(y)\leq f(z)$, and by the transitivity
of $\leq$ and by $(\ast)$ we obtain $xRz$. \qed

\vspace{0.3cm}

As an immediate consequence of Proposition 3 we have the following.

\vspace{0.3cm}

\noindent {\sc Proposition} 5. \quad {\it If $R,P\subseteq\W^2$ are preorders,
then $R=P$ iff $\F(R)=\F(P)$.}

\vspace{0.3cm}

\section{Representing \sp\ in \R}

Let \R\ be the category whose objects are sets in a certain universe,
and whose arrows are binary relations between these sets. Let
$I_a\subseteq a\times a$ be the identity relation on the set $a$, and
the composition $R_2\circ R_1\subseteq a\times c$ of $R_1\subseteq
a\times b$ and $R_2\subseteq b\times c$ is $\{(x,y)\mid(\exists z\in
c)(xR_1 z\;{\mbox{\it and }}zR_2 y)\}$.

Then let $p$ be a set in which the ordinal 2 is included, as in
the preceding section, and let the relation $\leq$ on $p$ be a
linear order such that every nonempty subset of $p$ has a greatest
element and 0 is the least element of $p$. A finite ordinal $p
\geq 2$ with the usual ordering satisfies these conditions, but in
general $p$ need not be finite ($p$ is up to isomorphism a
nonlimit ordinal with the inverse ordering). Then we define a map
$F_p$ from the objects of \sp\ to the objects of \R\ by setting
that $F_p(X)$ is $p^X$, namely, the set of all functions from $X$
to $p$. So $p^X$ has to be an object of \R.

For the functions $f_1:X\str p$ and $f_2:Y\str p$, let
$[f_1,f_2]:X^s\cup Y^t\str p$ be defined by
\[
[f_1,f_2](u)=\left\{
\begin{array}{lll}
f_1(s^{-1}_X(u)) & {\mbox{\rm if}} & u\in X^s
\\[.1cm]
f_2(t^{-1}_Y(u)) & {\mbox{\rm if}} & u\in Y^t.
\end{array}
\right.
\]

For $R:X\str Y$ an arrow of \sp, and for $f_1:X\str p$ and $f_2:Y\str p$
we define the arrow $F_p(R)$ of \R\ by
\[
(f_1,f_2)\in F_p(R)\quad {\mbox{\rm iff}}\quad [f_1,f_2]\in\F(R),
\]
where $\F(R)$ is the set of functions defined as in the preceding
section. Here \W\ is $X^s\cup Y^t$. Then we can prove the following
propositions.

\vspace{0.3cm}

\noindent {\sc Proposition} 6. \quad $F_p(\mj_X)=I_{F_p(X)}$.

\vspace{0.3cm}

\noindent {\it Proof}.\quad For $f_1,f_2:X\str p$ and $i,j\in\{s,t\}$ we have
$[f_1,f_2]\in\F(\mj_X)$ iff
\[
(\forall u\in X^i)(\forall v\in
X^j)(i^{-1}_X(u)=j^{-1}_X(v)\Rightarrow[f_1,f_2](u)\leq[f_1,f_2](v)).
\]
So if $[f_1,f_2]\in\F(\mj_X)$, then for $u=s_X(x)$ and $v=t_X(x)$ we have
$[f_1,f_2](u)\leq [f_1,f_2](v)$, which means $f_1(x)\leq f_2(x)$. By
setting $u=t_X(x)$ and $v=s_X(x)$ we obtain $f_2(x)\leq f_1(x)$, and by
the antisymmetry of $\leq$ we obtain $f_1(x)=f_2(x)$. So from
$[f_1,f_2]\in\F(\mj_X)$ we have inferred $f_1=f_2$.

\vspace{0.1cm}

For the converse, suppose $f_1=f_2$, and for $u\in X^i$ and $v\in X^j$
let $i^{-1}_X(u)=j^{-1}_X(v)$. If $i=j$, then $[f_1,f_2](u)\leq
[f_1,f_2](v)$ by the reflexivity of $\leq$. If $i\neq j$, then
$f_1(i^{-1}_X(u))=f_2(j^{-1}_X(v))$, and hence $[f_1,f_2](u)\leq
[f_1,f_2](v)$ by the reflexivity of $\leq$.
\qed

\vspace{0.3cm}

\noindent {\sc Proposition} 7.\quad $F_p(P\ast R)=F_p(P)\circ F_p(R)$.

\vspace{0.2cm}

\noindent {\it Proof}.\quad Suppose $R:X\str Y$ and $P:Y\str Z$.
We have to show that for $f_1:X\str p$ and $f_2:Z\str p$ such that
\[
(\ast\ast)\quad(\forall u,v\in X^s\cup Z^t)(u(P\ast R)v\Rightarrow
[f_1,f_2](u)\leq[f_1,f_2](v))
\]
there is an $f_3:Y\str p$ such that the following two statements are
satisfied:
\[
\begin{array}{l}
(\ast R)\quad(\forall u,v\in X^s\cup Y^t)(uRv\Rightarrow
[f_1,f_3](u)\leq[f_1,f_3](v)),
\\[.15cm]
(\ast P)\quad(\forall u,v\in Y^s\cup Z^t)(uPv\Rightarrow
[f_3,f_2](u)\leq[f_3,f_2](v)).
\end{array}
\]

For $y\in Y$ let
\[
\begin{array}{l}
X_y=_{\mbox{\scriptsize\it def}}
\{x\in X\mid (s_X(x),y)\in\tr(R^{-t}\cup P^{-s})\},
\\[.1cm]
Z_y=_{\mbox{\scriptsize\it def}}
\{z\in Z\mid (t_Z(z),y)\in\tr(R^{-t}\cup P^{-s})\}.
\end{array}
\]
Then we define $f_3$ as follows:
\[
f_3(y)=_{\mbox{\scriptsize\it def}}
\left\{\begin{array}{l}
\max(\{f_1(x)\mid x\in X_y\}\cup\{f_2(z)\mid z\in Z_y\})\;{\mbox{\rm if
}} X_y\cup Z_y\neq \emptyset
\\[.1cm]
0\quad{\mbox{\rm otherwise}}.
\end{array}
\right.
\]
We verify that $f_3$ satisfies $(\ast R)$. Suppose $uRv$.

\vspace{0.15cm}

(1) If $u,v\in X^s$, then $f_1(s^{-1}_X(u))\leq
f_1(s^{-1}_X(v))$ by $(\ast\ast)$, and hence $[f_1,f_3](u)\leq
[f_1,f_3](v)$.

\vspace{0.15cm}

(2) Suppose $u,v\in Y^t$. Since $uRv$, we must have
$X_{t^{-1}_Y(u)}\subseteq X_{t^{-1}_Y(v)}$ and
$Z_{t^{-1}_Y(u)}\subseteq Z_{t^{-1}_Y(v)}$, and hence
$f_3(t^{-1}_Y(u))\leq f_3(t^{-1}_Y(v))$. So $[f_1,f_3](u)\leq
[f_1,f_3](v)$.

\vspace{0.15cm}

(3) Suppose $u\in X^s$ and $v\in Y^t$. Since $uRv$, we must have
$s^{-1}_X(u)\in X_{t^{-1}_Y(v)}$, and hence $f_1(s^{-1}_X(u))\leq
f_3(t^{-1}_Y(v))$. So $[f_1,f_3](u)\leq [f_1,f_3](v)$.

\vspace{0.15cm}

(4) Suppose $u\in Y^t$ and $v\in X^s$.

\vspace{0.15cm}

(4.1) Let $X_{t^{-1}_Y(u)}\cup
Z_{t^{-1}_Y(u)}\neq \emptyset$ and $f_3(t^{-1}_Y(u))=f_1(x)$ for some
$x\in X_{t^{-1}_Y(u)}$. Since $uRv$, we must have
$(s_X(x),v)\in\tr(R^{-t}\cup P^{-s})$. Hence $(s_X(x),v)\in P\ast R$,
and we obtain $f_3(t^{-1}_Y(u))=f_1(x)\leq f_1(s^{-1}_X(v))$ by
$(\ast\ast)$. So $[f_1,f_3](u)\leq [f_1,f_3](v)$. We proceed
analogously when $f_3(t^{-1}_Y(u))=f_2(z)$ for some $z\in
Z_{t^{-1}_Y(u)}$.

\vspace{0.15cm}

(4.2) If $X_{t^{-1}_Y(u)}\cup Z_{t^{-1}_Y(u)}=\emptyset$,
then $f_3(t^{-1}_Y(u))=0$, and, since 0 is the least element of $p$, we
have $[f_1,f_3](u)\leq[f_1,f_3](v)$.

\vspace{0.15cm}

We verify analogously that $f_3$ satisfies $(\ast P)$.

\vspace{0.15cm}

It remains to show that if for some $f_3:Y\str p$ we have $(\ast R)$
and $(\ast P)$, then we have $(\ast\ast)$. Suppose $u(P\ast R)v$ and
$u,v\in X^s$. Then for $l\geq 0$ there is a (possibly empty) sequence
$y_1,\ldots,y_{2l}$ of elements of $Y$ such that
\[
(u,t_Y(y_1))\in R, (s_Y(y_1),s_Y(y_2))\in P,
(t_Y(y_2),t_Y(y_3))\in R,\ldots, (t_Y(y_{2l}),v)\in R;
\]
if $l=0$, then $uRv$.
Then by applying $(\ast R)$, $(\ast P)$ and the transitivity
of $\leq$ we obtain
$f_1(s^{-1}_X(u)) \leq f_1(s^{-1}_X(v))$, and hence
$[f_1,f_2](u) \leq [f_1,f_2](v)$. We proceed analogously when $u,v\in Z^t$,
or when $u\in X^s$ and $v\in Z^t$, or when $u\in Z^t$ and $v\in X^s$.
\qed

\vspace{0.3cm}

As an immediate consequence of Proposition 5 we have the following.

\vspace{0.3cm}

\noindent {\sc Proposition} 8.\quad {\it If $F_p(R)=F_p(P)$, then $R=P$.}

\vspace{0.3cm}

Since $F_p$ defined on the objects of \sp\ is clearly a one-one map,
Proposition 8 means that $F_p$ defined on the arrows of \sp\ is a
one-one map. So we have in \R\ as a subcategory an isomorphic
copy of \sp. Hence \sp\ is a category, as we promised we will show.
The two maps $F_p$, defined on the objects and on the arrows of
\sp, make a faithful functor from \sp\ to \R.

If the family \M\ satisfies the condition that for every $X,Y \in \M $
there is a $Z \in \M $ such that $Z$ is isomorphic to the
disjoint union of $X$ and $Y$, the category \sp\ has the structure of
a symmetric monoidal closed category (see \cite{ML71}, VII.7).

If we take the subcategory of \sp\ whose arrows correspond to binary
relations between members of \M, as explained at the end of Section 2,
and if \M\ is the set of objects of \R, then this subcategory of \sp\ is
isomorphic to \R. When our representation of \sp\ in \R\ via
the functor $F_p$ is restricted to this
subcategory, it amounts to a nontrivial embedding of \R\ in \R.

\section{Split preorders associated to proofs in\\ fragments of logic}

The language \L\ of conjunctive logic is built from a
nonempty set of propositional variables \P\ with the binary connective
$\wedge$ and the propositional constant, i.e. nullary connective, $\top$
(the exact cardinality of $\cal P$ is not important here). We use the
schematic letters $A,B,C,\ldots$ for formulae of \L.

We have the following axiomatic derivations for every $A$ and $B$ in
\L:
\[
\begin{array}{l}
\mj_A:A\str A,
\\
\k{\wedge}{1}_{A,B}:A\wedge B\str A,
\\
\k{\wedge}{2}_{A,B}:A\wedge B\str B,
\\
\k{\wedge}{}_A:A\str\top,
\end{array}
\]
and the following inference rules for generating derivations:
\[
\begin{array}{l}
\f{f:A\str B \quad g:B\str C}{g\circ f:A\str C}
\\[.3cm]
\f{f:C\str A \quad g:C\str B}{\langle f,g\rangle:C\str A\wedge B}
\end{array}
\]
This defines the deductive system \D\ of conjunctive logic (both
intuitionistic and classical). In this system $\top$ is included as an
``empty conjunction''.

Let \sp\ now be the category of split preorders of
$\omega=\{0,1,2,\ldots\}$. We define a map $G$ from \L\ to the objects
of \sp\ by taking that $G(A)$ is the number of occurrences of
propositional
variables in $A$. Next we define inductively a map, also denoted
by $G$, from the derivations of \D\ to the arrows of \sp:
\[
\begin{array}{ll}
{\makebox[2.5cm][l]{$G(\mj_A)=\mj_{G(A)},$}} &
\\[.1cm]
{\makebox[2.5cm][l]{$(u,v)\in G(\k{\wedge}{1}_{A,B})$}}  &
{\makebox[10cm][l]{${\mbox{\rm iff}} \quad
u=v\;{\mbox{\it or }}(u\in G(A\wedge B)^s\;{\mbox{\it and }}v\in
G(A)^t$}}
\\
 & {\makebox[10cm][l]{$\hspace{2.3cm}{\mbox{\it and }}s^{-1}_{G(A\wedge
B)}(u)=t^{-1}_{G(A)}(v)),$}}
\\[.1cm]
{\makebox[2.5cm][l]{$(u,v)\in G(\k{\wedge}{2}_{A,B})$}}  &
{\makebox[10cm][l]{${\mbox{\rm iff}} \quad
u=v\;{\mbox{\it or }}(u\in G(A\wedge B)^s\;{\mbox{\it and }}v\in
G(B)^t$}}
\\
 & {\makebox[10cm][l]{$\hspace{2.3cm}{\mbox{\it and }}s^{-1}_{G(A\wedge
B)}(u)=t^{-1}_{G(B)}(v)+G(A)), $}}
\\[.1cm]
{\makebox[2.5cm][l]{$(u,v)\in G(\k{\wedge}{}_{A})$}} &
{\makebox[10cm][l]{${\mbox{\rm iff}} \quad u=v,$}}
\end{array}
\]
\[
\begin{array}{ll}
{\makebox[2.5cm][l]{$G(g\circ f)=G(g)\ast G(f),$}} &
{\makebox[10cm][l]{}}
\end{array}
\]
\[
\begin{array}{ll}
{\makebox[2.5cm][l]{$(u,v)\in G(\langle f,g\rangle)$}}  &
{\makebox[10cm][l]{${\mbox{\rm iff}} \quad
u=v\;{\mbox{\it or }}(u\in G(C)^s\;{\mbox{\it and }}v\in
G(A\wedge B)^t$}}
\\
 & {\makebox[10cm][l]{$\hspace{2.3cm}{\mbox{\it and }}
 (u,t_{G(A)}(t^{-1}_{G(A\wedge B)}(v)))\in G(f))$}}
\\
 & {\makebox[10cm][l]{$\hspace{1.6cm}{\mbox{\it or }}
 (u\in G(C)^s\;{\mbox{\it and }}v\in G(A\wedge B)^t$}}
\\
 & {\makebox[10cm][l]{$\hspace{2.3cm}{\mbox{\it and }}
 (u,t_{G(B)}(t^{-1}_{G(A\wedge B)}(v)-G(A)))\in G(g)).$}}
\end{array}
\]

It is easy to check that for every arrow $f$ of \D, the relation $G(f)$ is a
preorder. When we draw this preorder, we may draw just the
corresponding strictification, as we remarked in Section 2. For
example, for $p,q,r\in \P$, the graph of $G(\langle \mj_{(p\wedge
q)\wedge \top}, \mj_{(p\wedge q)\wedge
\top}\rangle\circ\k{\wedge}{1}_{(p\wedge q)\wedge \top,r})$ would be

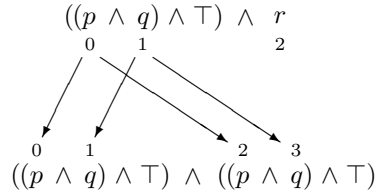
\begin{center}
\begin{picture}(120,80)

\put(17,55){\vector(-1,-2){15}}
\put(37,55){\vector(-1,-2){15}}
\put(24,55){\vector(3,-2){48}}
\put(44,55){\vector(3,-2){48}}

\put(10,10){\makebox(0,0){$(p\:\wedge\: q)$}}
\put(88,10){\makebox(0,0){$(p\:\wedge\: q)$}}
\put(30,70){\makebox(0,0){$(p\:\wedge\: q)$}}
\put(-8,10){\makebox(0,0){$($}}
\put(70,10){\makebox(0,0){$($}}
\put(12,70){\makebox(0,0){$($}}
\put(59,10){\makebox(0,0){$\wedge$}}
\put(40,10){\makebox(0,0){$\wedge\;\top)$}}
\put(118,10){\makebox(0,0){$\wedge\;\top)$}}
\put(60,70){\makebox(0,0){$\wedge\;\top)$}}
\put(79,70){\makebox(0,0){$\wedge$}}
\put(92,70){\makebox(0,0){$r$}}

\put(0,20){\makebox(0,0){\scriptsize$0$}}
\put(20,20){\makebox(0,0){\scriptsize$1$}}
\put(78,20){\makebox(0,0){\scriptsize$2$}}
\put(98,20){\makebox(0,0){\scriptsize$3$}}
\put(20,60){\makebox(0,0){\scriptsize$0$}}
\put(40,60){\makebox(0,0){\scriptsize$1$}}
\put(92,60){\makebox(0,0){\scriptsize$2$}}

\end{picture}
\end{center}
where the formulae of \L\ are written down to show where the ordinal 3
of the source and the ordinal 4 of the target come from. Actually, we
may as well omit the ordinals from such drawings.

We may also replace uniformly the relations defined above
as values of $G$ by the converse
relations. This shows that the information about the orientation of the
edges is not essential. We may as well omit this information when we
draw graphs, provided we take care that in composing graphs edges of
a single graph are not composed with each other.

If we stipulate that for $f,g:A\str B$ derivations of \D\ we have that
$f$ and $g$ are equivalent iff $G(f)=G(g)$, and define a {\it proof}
of \D\ to be the equivalence class of a derivation, then proofs of \D\
would make the arrows of a category \C, whose objects are formulae of
\D, with the obvious sources and targets. In this particular case,
where \D\ is our deductive system for conjunctive logic, the
category \C\ will be the free cartesian category generated by the set
of propositional variables \P\ as the generating set of objects (this
set may be conceived as a discrete category). This fact about \C\
follows from the coherence result for cartesian categories treated in
\cite{DP01} and \cite{P02}. (Cartesian categories are categories with
all finite products, including the empty product, i.e. terminal
object. The category \C\ can be equationally presented; see
\cite{LS86}, Chapter I.3, or \cite{DP01}.)

When we replace the split
preorders $G(f)$ defined above by their transitive and symmetric
closures, we obtain the equivalence relations of \cite{DP02a}, which
are also split preorders. The category \C\ induced by these split
preorders is, however, again the free cartesian category generated by
\P.

The language of disjunctive logic is dual to the language we had above:
instead of $\wedge$ and $\top$
we have $\vee$ and $\bot$ in $\cal L$, and instead
of $\k{\wedge}{\it i}$, $\k{\wedge}{}$ and $\langle\; ,\; \rangle$ we have
in the corresponding deductive system $\cal D$
\[
\begin{array}{l}
\k{\vee}{1}_{A,B}:A\str A\vee B,
\\
\k{\vee}{2}_{A,B}:B\str A\vee B,
\\
\k{\vee}{}_A:\bot\str A,
\end{array}
\]
\[
\f{f:A\str C \quad g:B\str C}{[ f,g]:A\vee B\str C}
\]
We obtain that equivalence of derivations induced by maps $G$ dual to
those above makes of the proofs of \D\ the free category with all
finite coproducts generated by \P.

Let us now assume we have in \D\ both $\wedge$ and $\vee$, but without
$\top$ and $\bot$, and let us introduce the category \C\ of proofs of
\D\ as we did above. To define the new map $G$,
we combine the two kinds of map $G$ defined
previously, paying attention to order (for example, both in
$G(\k{\wedge}{\it i}_{A,B})$ and $G(\k{\vee}{\it i}_{A,B})$
edges connecting the domain and the codomain
must be oriented in the same direction).
The category \C\ will be the free category with nonempty
finite products and coproducts generated by \P. This follows from the
coherence result of \cite{DP02}. In the presence of $\top$ and $\bot$,
we would obtain a particular brand of bicartesian category, according
to the coherence result of \cite{DP01a}. If we replace the split
preorders $G(f)$ considered here for conjunctive-disjunctive logic by
their transitive and symmetric closures, as we did above for
conjunctive logic, we will not obtain any more as the category \C\ the
free category with nonempty finite products and coproducts (see
\cite{DP02a}).

Up to now, in this section, the split preorders that were values
of $G$ always corresponded to relations between finite ordinals
(see the end of Sections 2 and 4). We will next consider
split preorders that are not such.

The coherence result of \cite{KML71} for symmetric monoidal closed
categories without the unit object $I$ can be used to show that split
preorders for the appropriate deductive system, which corresponds to
the tensor-implication fragment of intuitionistic linear logic, give
rise to the free symmetric monoidal closed category without $I$. These
split preorders may be taken as equivalence relations, but they need
not be such.

Suppose we add intuitionistic implication $\rightarrow$ to conjunctive
or conjunctive-disjunctive logic, and take the equivalence of
derivations induced by split preorders that for the derivations from
$p\wedge(p\rightarrow q)$ to $q$ and from $q$ to $p\rightarrow(p\wedge
q)$ would correspond to the following graphs

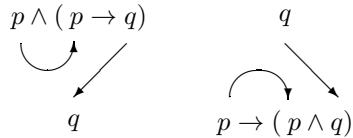
\begin{center}
\begin{picture}(120,60)

\put(40,40){\vector(-1,-1){20}}

\put(20,39){\vector(0,1){2}}

\put(10,40){\oval(20,20)[b]}

\put(20,10){\makebox(0,0){$q$}}
\put(22,50){\makebox(0,0){$p\wedge(\:p\rightarrow q)$}}

\put(100,40){\vector(1,-1){20}}

\put(101,21){\vector(0,-1){2}}

\put(90,20){\oval(22,20)[t]}

\put(100,50){\makebox(0,0){$q$}}
\put(100,10){\makebox(0,0){$p\rightarrow(\:p\wedge q)$}}

\end{picture}
\end{center}

\noindent In this case, however, the category \C\ would not be the free
cartesian closed or free bicartesian closed category
generated by \P\ (for counterexamples see \cite{DP02}, Section 1, and
\cite{S75}). Taking the transitive and symmetric closures of our split
preorders would again not give this free cartesian closed or free
bicartesian closed category.

In \cite{DP02a} we considered the category \G\ of split preorders on
$\omega$ that are equivalence relations, and represented \G\ in a
subcategory of \R\ by a functor that amounts to a particular case of
the functor $F_p$ considered in this paper. In \cite{DP02a}, the set $p$
is a finite ordinal, and $F_p(n)=p^n$ is taken to be a finite ordinal too.
This representation of \G\ is connected to Brauer's representation of
Brauer algebras of \cite{B37}, as explained in \cite{DP02a} (Section
6). This is the reason why we call our representation of split preorders
{\it Brauerian}. It is a generalization of Brauer's representation, and
of the Brauerian representation of \G\ of \cite{DP02a}.

\vspace{.5cm}

\noindent {\footnotesize {\it Acknowledgement.} The writing of this paper
was financed by the Ministry of Science, Technology and Development of
Serbia through grant 1630 (Representation of proofs with applications,
classification of structures and infinite combinatorics).}

\end{document}